\documentclass[12pt]{article}

\newtheorem{Le}{Lemma}
\usepackage[latin1]{inputenc}
\newtheorem{Co}{Conjecture}

\begin{document}

\title{Erd\"{o}s's Matching Conjecture and $s$-wise $t$-intersection Conjecture}

\date{July 10,  2013}
\author{Vladimir Blinovsky\footnote{The author was supported by 
  NUMEC/USP (Project MaCLinC/USP).
}}
\date{\small
 Instituto de Matematica e Estatistica, USP,\\
Rua do Matao 1010, 05508- 090, Sao Paulo, Brazil\\
 Institute for Information Transmission Problems, \\
 B. Karetnyi 19, Moscow, Russia,\\
vblinovs@yandex.ru}

\maketitle\bigskip

\begin{center}
{\bf Abstract}
\bigskip

We find the formula for the maximal cardinality of the family of $n$-tuples from ${[n]\choose k}$ with does not have $\ell$--matching. This  formula after some analytical issues can be reduce to the  Erd\"{o}s's Matching formula. Also we prove the conjecture about the cardinality of maximal $s$-wise $t$-intersecting family of $k$-element subsets of $[n]$. In the proofs we use original method which we have already used  in the proof of Mikl\'{o}s-Manikam-Singhi conjecture in~\cite{1}.
 \end{center}

\bigskip

{\bf I Introduction and Formulation of Results}
\bigskip

Define $[n]=\{ 1,\ldots ,n\}$ and ${[n]\choose k}=\{ E\subset [n]:\ |E|=k\}$. We say that family ${\cal A}\subset {[n]\choose k}$ has $\ell$-matching if there exists the set $\{ E_i ,\ i\in [\ell ]\}\subset {\cal A}$ such that $E_i \bigcap E_j =\emptyset$ when $i\neq j$. 

First problem which we would like to introduce is to find the maximal cardinality $M(\ell ,n ,k)$ of ${\cal A}\subset {[n]\choose k}$ which has no $\ell$-matching.

In 1965 Erd\"{o}s~\cite{2} formulate the following 
\begin{Co}
\label{co1}
The value $M(\ell , n,k)$ satisfies the following equality
\begin{equation}
\label{e1}
M(\ell ,n,k)=\max\left\{ {k\ell -1\choose k},{n\choose k}-{n-\ell +1\choose k}\right\} .
\end{equation}
\end{Co}
This conjecture is one of the main statements in extremal hypergraph theory.
Erd\"{o}s wrote in~\cite{2} that he manage to prove this corollary for $k=2$ and for $\ell =2$ it is Erd\"{o}s-Ko-Rado
result but general case seems elusive. 

Later this corollary was confirmed for several conditions on parameters of the problem, we mention the 
proof of the conjecture for $n\geq(2\ell -1)k-\ell+1$ in~\cite{4}, there was proved that in this case
$$
M(\ell ,n,k)={n\choose k}-{n-\ell +1\choose k} .
$$
Also the conjecture was proved for $k=3$ in~\cite{5}.

Let's mention also that asymptotic equality for $M(\ell ,n,k)$ , which follows from the conjecture is proved for some parameters in~\cite{6}.  

First our result is the proof of the following 

\begin{Le}
\label{le2}
The following equality is valid
\begin{equation}
\label{eq1}
M(\ell ,n,k)=\max_{1\leq i\leq k}\sum_{j\geq i}{\ell i-1\choose j}{n-\ell i+1\choose k-j} .
\end{equation}
\end{Le}

Thus the proof that the conjecture~\ref{co1} is true for  all parameters $\ell ,n,k$ reduced to the proof 
of technical equality
\begin{eqnarray*}
&&
\max_{1\leq i\leq k}\sum_{j\geq i}{\ell i-1\choose j}{n-\ell i+1\choose k-j}\\
&=&
\max\left\{ {k\ell -1\choose k},{n\choose k}-{n-\ell +1\choose k}\right\} .
\end{eqnarray*}

Note, that for the arbitrary $i\in [k]$  the choice of the set 
$$
{\cal A}=\left\{ A\in{[n]\choose k}:\ |A\bigcap [\ell i-1]|\geq i \right\}
$$
shows
that
$$
M(\ell ,n,k)\geq\max_{1\leq i\leq k}\sum_{j\geq i}{\ell i-1\choose j}{n-\ell i+1\choose k-j} .
$$
So we need to prove the opposite inequality

To introduce our second result we introduce some additional notations. We say that family ${\cal B}\subset {[n]\choose k}$ is $s$-wise $t$-intersecting if for the arbitrary subset $\{ E_i , i\in [s]\}\subset{\cal B}$ the following relation is true $|E_1 \bigcap E_2 \bigcap\ldots  \bigcap E_s  |\geq t$. 
Let $N(s,n,k,t)$ is the maximal cardinality of $s$-wise $t$-intersecting family from ${[n]\choose k}$.

There is the following old 
\begin{Co}
\label{co2} Let $sk<(s-1)n +t$.The following equality is valid
$$
N(s,n,k,t)=\max_{r\geq 0}\left\{ \biggl| E\in {[n]\choose k}:\ |E\bigcap [t+rs]|\geq t+(s-1)r\biggr | \right\} .
$$
\end{Co}
Note that choice of the set
$$
\left\{ \biggl| E\in {[n]\choose k}:\ |E\bigcap [t+rs]|\geq t+(s-1)r\biggr | \right\}
$$
for $r\geq 0$  shows that
$$
N(s,n,k,t)\geq\max_{r\geq 0}\left\{ \biggl| E\in {[n]\choose k}:\ |E\bigcap [t+rs]|\geq t+(s-1)r\biggr | \right\} .    
$$
So we need to prove the opposite inequality.

Note also that if $s k\geq (s -1)n +t$, then the whole set ${[n]\choose k}$ is $s$ - wise $t$ - intersecting family.
There are many publications which are devoted to solution of this problem in particular cases.
The most important result was obtained by Ahlswede and Khachatrian in celebrated paper~\cite{7}. They confirm the validness of this conjecture for the case $s=2$. In all other cases there are partial solutions (when some parameters are given and $n$ is sufficiently large).
We mention papers~\cite{9}-\cite{10}. 

Our second result is the proof of this conjecture for all parameters $s ,n,k$.

We note that in~\cite{11} we prove the fractional analog of lemma~\ref{le2}. Hence we have confirmed the expression
similar to~(\ref{eq1}) for fractional matching also. 

The paper is organized as follows: in Section~II we introduce the {\it Symmetrical Smoothing Method}, which actually we have already used in~\cite{1} and in~\cite{11}. We also formulate and prove technical lemma~\ref{le2} which we used later in section III. In Section III we, using lemma~\ref{le2}, complete the proof of Conjecture~\ref{co2}. \bigskip

{\bf Section II}
\bigskip

Next we use the natural bijection between $2^{[n]}$ and set of binary $n$-tuples $ \{ 0,1\}^n$ and make no difference between these two sets.

We say that family ${\cal A}\subset {[n]\choose k}$ is (left) compressed if from the inclusion $A=(a_1 ,\ldots , a_k  )\in {\cal A}$ and the conditions $b_j \leq a_j $ it  follows that $B=(b_1 ,\ldots ,b_k )\in{\cal A} $.
Note, that we can assume that the extremal intersection families are left compressed. 

Also we assume that left compressed family ${\cal A}\subset{[n]\choose k}$ is defined by the inequalities
\begin{equation}
\label{k1}
{\cal A}=\left\{ x\in{[n]\choose k}:\ (\omega_i ,x) > 0,\ i\in [N]\right\} ,
\end{equation}
where
$$\omega_i =(\omega_{i,1},\ldots ,\omega_{i,n})\in R^n
$$
and 
\begin{equation}
\label{g1}
\omega_{i,j}\geq\omega_{i,j+1},\ j\in [n-1],\ \sum_{j=1}^n \omega_{i,j}=1.
\end{equation}
Indeed, the set ${\cal A}$ is shifted. Arbitrary left compressed set can be defined as ${\cal A}$. However as we will see later in our case we can restrict ourselves assuming that the 
extremal set can be generated by one inequality only. 

Next we assume that extremal families in both problems are left compressed and are defined by one inequality.  It is easy to see that if the family ${\cal A}$ has $\ell$ - matching, then the non intersecting set $(x_1 ,\ldots ,x_\ell)\subset {\cal A}$ can be chosen in such a way that $x_i \subset [\ell k]$.

The Symmetrical Smoothing Method consists in approximation of the number $|{\cal A}|$ by the smooth symmetric function $\omega$, which allows to use analytic methods to determine the values $\omega_j$ on which achieves extremum of $|{\cal A}|$.

Some of the values $\omega_j$ can be negative. Next we make transformations of $\omega$ and write the inequality in the determination~(\ref{k1}) of ${\cal A}$ in the equivalent form, where all coefficients are nonnegative. Consider the following set of basis vectors $z_j =(k\ell -dj,\ldots ,k\ell-dj,-dj,-dj,\ldots ,-dj)\in R^n$, where the number of coordinates $k\ell-dj$ is equal to $j$ and $j\in [k\ell-1]$. Because the maximal set ${\cal A}$ is compressed, we need to choose only first $k\ell$ coordinates $\omega_{j}$ and other coordinates we can choose as large as possible, i.e. all of them are equal to $\omega_{k\ell}$.

Then it is easy to check, that vector $\omega$, which coordinates satisfies inequalities~(\ref{g1}) and which determine the maximal family in the first problem, can be represent as the sum
\begin{equation}
\label{gh2}
\omega =\sum_{j=1}^{k\ell -1}\alpha_j z_j
\end{equation}
with nonnegative coordinates $\alpha_j \geq 0$ and some $d$. Indeed from~(\ref{gh2}) follows that for $j\leq k\ell -1$
$$
\omega_j - \omega_{j+1}=\alpha_j k \ell
$$
or
$$
\alpha_j =\frac{\omega_j -\omega_{j+1}}{k\ell}\geq 0.
$$
The last equation contains only difference of $\omega_j$ we have degree of freedom to determine $\omega_{k\ell}$. To make this we choose proper $d$. 
We have
$$
\omega_i =-d\sum_{j=1}^{k\ell -1}j\alpha_j +k\ell\sum_{j=i}^{k\ell-1}\alpha_j
$$
and
$$
\lambda = \sum_{i=1}^{k\ell}\omega_i =-d\frac{k\ell-1}{k\ell}(\lambda -k\ell\omega_{k\ell}) +\lambda - k\ell \omega_{k\ell}
$$
and hence
$$
d=-\frac{k^2\ell^2\omega_{k\ell}}{(k\ell -1)(\lambda -k\ell\omega_{k\ell})}.
$$

Substituting expansion~(\ref{gh2}) to the inequality $(\omega ,x)>0$, we obtain
\begin{eqnarray}
\label{ds1}
&& \left(\sum_{j=1}^{k\ell -1}\alpha_j z_j ,x\right) =\sum_{j=1}^{k\ell-1}\alpha_{j}\left(k\ell \sum_{m=1}^j x_m -jkd\right)\\
\nonumber && =k\ell\sum_{j=1}^{k\ell -1}\left(\sum_{m=j}^{k\ell -1}\alpha_m \right) x_j -kd\sum_{j=1}^{k\ell -1}j\alpha_j  >0,\ \delta =\frac{d}{\ell}= -\frac{k^2\ell\omega_{k\ell}}{(k\ell -1)(\lambda -k\ell\omega_{k\ell})}\end{eqnarray}

From inequality $\delta >0$ it follows that $\omega_{k\ell} <0$. Because inequalities~(\ref{k1}) which determine set ${\cal A}$ are strict, varying coordinates of $\omega$ we can 
we can reach the case, when $\delta$ is rational number and inequalities in~(\ref{k1}) with new rational $\delta$ are still determine ${\cal A}$. 

Define
$$
\beta_j =\frac{\sum_{m=j}^{k\ell-1}\alpha_m}{\sum_{m=1}^{k\ell -1}m\alpha_m} =\frac{\omega_j - \omega_{k\ell}}{\lambda-k\ell\omega_{k\ell}},\ j\in [k\ell-1].
$$
We have
$$ 
\beta_{i,j}\geq 0,\ \beta_{i,j}=0,\ j> k\ell-1;\ \sum_{j=1}^n \beta_{i,j} =1,\ \beta_{i,j}\geq \beta_{i,j+1}.
$$
We rewrite definition of ${\cal A}$ as follows
\begin{equation}
\label{dd1}
{\cal A}=\left\{ x\in {[n]\choose k}:\ \langle\beta_i ,x\rangle >\delta_i (x),\ i\in [N]\right\},
\end{equation}
where  $\delta_i (x)$ depends on $x$.

Define function
$$
\varphi (\{\beta \},x)\stackrel{\Delta}{=}\frac{1}{(2\pi )^{N/2}}\prod_{i=1}^N \int_{-\infty}^{(\langle\beta_i ,x\rangle -\delta_i (x))/\sigma_i (m)}e^{-\xi^2/2}d\xi\to\left\{\begin{array}{ll}
1,& x\in{\cal A};\\
0,& \hbox{otherwise}
\end{array}
\right.
$$
as $\sigma_i (m)\to 0$ and $\delta_i (x)\to \delta_j$ as $m\to \infty$. 

We have
$$
\Biggl| |{\cal A}|-\sum_{x\in {[n]\choose2}}\frac{1}{(2\pi)^{N/2}}\prod_{i=1}^N \int_{-\infty}^{(\langle\beta_i ,x\rangle -\delta_i (x))/\sigma_i (m)}e^{-\xi^2/2}d\xi \Biggr|\to 0
$$
as $m\to\infty$. 

Denote
$$
A(x,\{\omega\},i)=\frac{1}{(2\pi )^{N/2}}\prod_{j=1,\ j\neq i}^N \int_{-\infty}^{(\langle x,\beta_j\rangle -\delta_j (x))/\sigma_j (m)}e^{-\xi^2 /2}d\xi.
$$
We have
\begin{equation}
\label{gg1}
\varphi (\{\beta\},x) \stackrel{\Delta}{=} A(x,\{\omega\},i)\int_{-\infty}^{(\langle x,\beta_i\rangle -\delta_i (x))/\sigma_i (m)}e^{-\xi^2/2}d\xi .
 \end{equation}
 We have partial derivative
 \begin{eqnarray}
 \label{kk2}&&
 \left(\sum_{x\in {[n]\choose2}}\frac{1}{(2\pi)^{N/2}}\prod_{i=1}^N \int_{-\infty}^{(\langle\beta_i ,x\rangle -\delta_i (x))/\sigma_i (m)}e^{-\xi^2/2}d\xi \right)^\prime_{\beta_{i,j}}\\&&\nonumber
  =\frac{1}{\sigma_i(m)}\Biggl( \sum_{x\in {[n]\choose k}:\ x_j=1,\ x_a =0}A(x,\{\omega\},j)e^{-\frac{(\langle\beta_i ,x\rangle -\delta_i (x))^2}{2\sigma_i^2 (m)}}\\ &&
 -\sum_{x\in {[n]\choose k}:\ x_j=0,\ x_a =1}A(x,\{\omega\},i)e^{-\frac{(\langle\beta_i ,x\rangle -\delta_i (x))^2}{2\sigma_i^2 (m)}}\Biggr) .\nonumber
 \end{eqnarray}
  Denote
 \begin{eqnarray*}
 &&
 x(r)=\hbox{arg}\min_{x\in {[n]\choose k}\setminus\left\{\bigcup_{i=1}^{r-1}x(i)\right\},\ \{ x_j =1,\ x_a =0\},\ \{ x_j =0,\ x_a =1\}} (\langle \beta_j ,x\rangle -\delta_j (x))^2;\\
 &&  
\bar{x}(r)=\hbox{arg}\min_{x\in {[n]\choose k}\setminus\left\{\bigcup_{i=1}^{r-1}x(i)\right\},\ \{ x_j =0,\ x_a =1\},\ \{ x_j =1,\ x_a =0\}} (\langle \beta_j ,x\rangle -\delta_j (x))^2 .\end{eqnarray*}  
We have
$$
A(x(r),\{\beta\},i)=e^{-\frac{\ln (A(x(r),\{\beta\},i)\sigma_j^2 (m)}{\sigma_i^2 (m)}},\ A(\bar{x}(r),\{\beta\},i)=e^{-\frac{\ln (A(\bar{x}(r),\{\beta\},i))\sigma_j^2 (m)}{\sigma_i^2 (m)}}
$$
 and $\ln (A(x(r),\{\beta\},i))\sigma_i^2 (m),\ \ln(A(\bar{x}(r),\{\beta\},i))\sigma_i^2 (m)\to 0$ as $m\to\infty$. 
 
 We can choose $\delta_i (x(r), \sigma_i (m)),\ \delta_i (\bar{x}(r),\sigma_i (m))$ s.t.
 \begin{eqnarray}
 \label{l1}&&
 A(x(r),\{\beta\},i)e^{-\frac{(\langle\beta_i ,x(r)\rangle - \delta_i (x(r)))^2}{2\sigma_i^2 (m)}} = e^{-\frac{(\langle\beta_i ,x(r)\rangle - \delta_i (x(r), \sigma_i (m)))^2}{2\sigma_i^2 (m)}};\\
 && 
 A(\bar{x}(r),\{\beta \},i)e^{-\frac{(\langle\beta_i ,\bar{x}(r)\rangle - \delta_i (\bar{x}(r)))^2}{2\sigma_i^2 (m)}} = e^{-\frac{(\langle\beta_i ,\bar{x}(r)\rangle - \delta_i (\bar{x}(r), \sigma_i(m)))^2}{2\sigma_i^2 (m)}}\label{l2}
 \end{eqnarray}
 and
$$
 e^{-\frac{(\langle\beta_i ,x(r)\rangle - \delta_i (x(r), \sigma_j(m)))^2}{2\sigma_i^2 (m)}}=e^{-\frac{(\langle\beta, \bar{x}(r)\rangle - \delta (\bar{x}(r), \sigma_i (m)))^2}{2\sigma^2 (m)}}
 $$
 and
  $$
  \lim_{m\to\infty}\delta (x,\sigma (m))=\delta.
  $$
We have
 \begin{equation}
 \label{kk1}
 (\langle\beta ,x(r)\rangle -\delta)^2 = (\langle\beta ,\bar{x}(r)\rangle -\delta)^2. 
 \end{equation}
 Summing last inequalities over $r$ we obtain equality
 \begin{eqnarray}&&
 \label{bb1}
 \sum_{x\in {[n]\choose k}:\ x_i=1,\ x_a=0}  \left(\beta_{i}  +\sum_{p=1,p\neq i}^{a-1}\beta_{p} x_p -\delta\right)^2 \\
 &&
  =\sum_{x\in {[n]\choose k}:\ x_i=0,\ x_a=1}  \left(\beta_{a}  +\sum_{p=1,p\neq i}^{a-1}\beta_{p} x_p -\delta\right)^2 \nonumber
  \end{eqnarray}
  or
  $$
  (\beta_{i}-\beta_{a})\left((\beta_{i}+\beta_{a}-2\delta ){n-2\choose k-1}+2{n-3\choose k-2}\sum_{p=1,p\neq i}^{a-1}\beta_{p}\right)=0
  $$
  or
  $$
  \beta_{i}=\beta_{a},$$
    \begin{eqnarray}&&\label{re1}
  \beta_{i}+\beta_{a}= \frac{\delta (n-2)-k+1}{\frac{n}{2} -k}=2\delta -2\frac{(1-2\delta)(k-1)}{n-2k}.
  \end{eqnarray}
Assume that $\beta_{i}=\alpha,\ i\in [m],\ \beta_{i}=\beta,\ i\in [m+1,a],\ \beta_{i}=0,\ i>a$.   Then using~(\ref{re1}), we have system of equations
  \begin{eqnarray}
  \label{a1}&&
  \alpha m+\beta (a-m)=1,\\
  && \alpha +\beta = 2\delta -2\frac{(1-2\delta)(k-1)}{n-2k} .
 \label{aa2} 
 \end{eqnarray}
  Assume that $\beta_{i}=\alpha,\ j\in [m],\ \beta_{i}=\beta,\ i\in [m+1,a],\ \beta_{i}=0,\ i>a$.  assume also that from at least one of equalities~(\ref{kk1}) follows equality
  \begin{equation}\label{ff1}
   \langle \beta,x\rangle +\langle \beta,\bar{x}\rangle=2\delta .
   \end{equation}
Hence there exist vectors $x ,\bar{ x}  \in {[n]\choose k}$, s.t. 
 \begin{eqnarray*}&& x_{i}=1,\ x_{a}=0,\ \bar{x}_i=0, \bar{x}_a =1,\\&&
 |\{x_{1},\ldots x_{m}\}| |= q_1,\   |\{x_{m+1},\ldots x_{a}\}| =q_2 ,\\&& |\{ \bar{x}_{1},\ldots ,\bar{x}_{m}\}| = \bar{q}_1,  |\{ \bar{x}_{m+1},\ldots ,\bar{x}_{a}\}| = \bar{q}_2 .
 \end{eqnarray*}
   We can rewrite equality~(\ref{ff1}) as follows 
   $$
   \alpha (q_1 +\bar{q}_1)+\beta (q_2+\bar{q}_2)=2\delta.
   $$
   Denote
   $$
   w=q_1+\bar{q}_1,\ \bar{w}=q_2+\bar{q}_2.
   $$
   Using equality~(\ref{aa2}) we have
   $$
   \alpha\left(2i-\frac{n-2k}{n-2}\right)+\beta \left(2j-\frac{n-2k}{n-2}\right) =\frac{2(k-1)}{n-2}.
   $$
   From inequalities 
   $$
   2i-\frac{n-2k}{n-2} >i,\ 2j-\frac{n-2k}{n-2} >j.   $$
   we have
   $\delta <\frac{k-1}{n-2}$. From other side from~(\ref{aa2}) follows inequality
   $$
   \delta>\frac{k-1}{n-2}.
   $$
   This contradiction allows us to skip possibility of the case~(\ref{ff1}). 
   
    Hence we can assume that from~(\ref{kk1}) follows equalities
   $$
   \langle \beta ,x(r)\rangle = \langle \beta ,\bar{x}(r)\rangle .
   $$
   Summing these inequalities over $r$, we obtain relation
   $$\sum_{r=1}^{n-2\choose k-1}\langle \beta_i ,x(r)\rangle = \sum_{r=1}^{n-2\choose k-1}\langle \beta_i ,\bar{x}(r)\rangle 
   $$
   from which easily follows equality
   $$\alpha =\beta .$$
   Hence it is sufficient to consider the case $\alpha=\beta =\frac{1}{a}.$ 
    
We can rewrite inequalities in the definition of ${\cal A}$  in~(\ref{dd1}) as follows.

Let's 
$$
\beta_{i}\geq 0,\ \beta_i =0,\ i >k\ell-1;\ \sum_{i=1}^n \beta_i =1,\ \beta_{i}\geq \beta_{i+1}.
$$
Then the set
\begin{equation}
\label{e09}\left\{ x\in{[n]\choose k}:\ (\beta ,x) >\delta \right\}
\end{equation}
is compressed. 
It is easy to see, that arbitrary compressed set ${\cal A}$ is the intersection of finite number $N$ of compressed sets.

To prove  lemma~\ref{le2}  we formulate the basic Optimization problem~1.
\bigskip

{\bf Optimization problem 1.}
\bigskip

Find maximum over the choice of one-step functions $\{ \beta_{i}\}$ and $\{\delta_i \}$ of the function
$$
|{\cal A} |
 $$
 under the condition that some special choice of matching $\{ x_s \}$ (which we choose later) does not belong to ${\cal A}$.
 
  We choose matching set $X =\{ x_1 ,\ldots ,x_\ell \}$ of nonintersecting elements from ${[n]\choose k}$:
 $$
 x_j =\{ j , j+\ell ,\ldots , j+(k-1)\ell \} ,\ j\in [\ell ].
 $$
 This set is forbidden to be included in ${\cal A}$ by  one-step function $\beta$, i.e. for some  $j\in [\ell ]$
 \begin{equation}
 \label{es1}
 (\beta ,x_j )> \delta  .
 \end{equation}
 Let's  $a = (m-1)\ell +p$, for some $m\in [k]$ and $p \in [0,\ell )$ and $\beta_{q}=1/a $ when $q\in [a]$. The restriction~(\ref{es1}) means that it is  sufficient to choose $\delta<\psi$ where
\begin{equation}
\label{et11}
\psi =\frac{m-1}{a}=\frac{m-1}{(m-1)\ell+p}.
\end{equation}
For this choice $\beta$ and $\delta$, which is sufficiently close to $\psi$ we have 
\begin{equation}
\label{e4r}
\sum_{j\geq m}{(m-1)\ell+p\choose j}{n-(m-1)\ell -p\choose k-j}
\end{equation}
choices of admissible $x\in{[n]\choose k}$. 
It is left to find the maximum over choices of values of $a \in [k\ell -1]$ of the sum from~(\ref{e4r}):
\begin{eqnarray}
\label{edd1}
&&\max_{a\in [k\ell -1]}\sum_{j> a/\ell}{a\choose j}{n-a\choose k-j}\\
&=& \nonumber
\max_{i\in [k]}\sum_{j\geq i}{\ell i -1\choose j}{n-\ell i +1\choose k-j} .
 \end{eqnarray}  
Equality in~(\ref{edd1}) follows from the fact, that 
$$
\sum_{j >a/\ell }{a\choose j}{n-a\choose k-j}
$$
decreases as $a$ decreases from $\ell i-1$ to $\ell (i-1)$.

Lemma~\ref{le2} is proved.
\bigskip
\newpage

{\bf Section III. Proof of the Conjecture~\ref{co2}}
\bigskip

Next we make transformations of $\omega$ and write the inequality in the determination~(\ref{k1}) of ${\cal A}$ in the equivalent form, where all coefficients are nonnegative. Consider the following set of basis vectors $z_j =(n -dj,\ldots ,n-dj,-dj,-dj,\ldots ,-dj)\in R^n$, where the number of coordinates $n-dj$ is equal to $j$ and $j\in [n-1]$. Because the maximal set ${\cal A}$ is compressed, we need to choose first $n-1$ coordinates $\omega_{j}$ and $\omega_n$ we can choose as large as possible, i.e.  $\omega_{n-1}=\omega_{n}$.

Then it is easy to check, that vector $\omega$, which coordinates satisfies inequalities~(\ref{g1}) and which determine the maximal family in the first problem, can be represent as the sum
\begin{equation}
\label{gh2}
\omega =\sum_{j=1}^{n -1}\alpha_j z_j
\end{equation}
with nonnegative coordinates $\alpha_j \geq 0$ and some $d$. Indeed from~(\ref{gh2}) follows that for $j\leq n -1$
$$
\omega_j - \omega_{j+1}=\alpha_j n
$$
or
$$
\alpha_j =\frac{\omega_j -\omega_{j+1}}{n}\geq 0.
$$
We have
$$
\omega_i = -d\sum_{j=1}^{n-1} j\alpha_j +n\sum_{j=i}^{n-1}\alpha_j ,\ i\leq n-1
$$
and hence
$$
\lambda\stackrel{\Delta}{=}\sum_{i=1}^{n-1}\omega_i =-d\frac{n-1}{n}(\lambda -n\omega_n)+\lambda-n\omega_n 
$$
or
$$
d=-\frac{n^2}{n-1}\frac{\omega_n}{\lambda-n\omega_n}.
$$
Substituting expansion~(\ref{gh2}) to the inequality $(\omega ,x)>0$, we obtain
\begin{eqnarray}
\label{ds1}
&& \left(\sum_{j=1}^{n -1}\alpha_j z_j ,x\right) =\sum_{j=1}^{n-1}\alpha_{j}\left(n \sum_{m=1}^j x_m -jkd\right)\\
\nonumber && n\sum_{j=1}^{n -1}\left(\sum_{m=j}^{n -1}\alpha_m \right) x_j -kd\sum_{j=1}^{n -1}j\alpha_j  >0,\ \delta =\frac{kd}{n}=-\frac{kn}{n-1}\frac{\omega_n}{\lambda-n\omega_n}.\end{eqnarray}
From inequality $\delta >0$ it follows that $\omega_{n} <0$. Because inequalities~(\ref{k1}) which determine set ${\cal A}$ are strict, varying coordinates of $\omega$ we can 
we can reach the case, when $\delta$ is rational number and inequalities in~(\ref{k1}) with new rational $\delta$ are still determine ${\cal A}$. 

Define
$$
\beta_i=\frac{\sum_{m=i}^{n-1}\alpha_m}{\sum_{m=1}^{n -1}m\alpha_m} =\frac{\omega_i - \omega_{n}}{1-n\omega_{n}},\ i\in [n-1].
$$
We have
$$ 
\beta_{i}\geq 0,\ \beta_{n}=0;\ \sum_{j=1}^n \beta_{i} =1,\ \beta_{i}\geq \beta_{i+1}.
$$

To prove conjecture~\ref{co2} we will follow the same procedure as in the previous section.

$$
  \beta_{i}=\beta_{a},$$
    \begin{eqnarray}&&\label{re1}
  \beta_{i}+\beta_{a}= \frac{\delta (n-2)-k+1}{\frac{n}{2} -k}=2\delta -2\frac{(1-2\delta)(k-1)}{n-2k}.
  \end{eqnarray}
Assume that $\beta_{i}=\alpha,\ i\in [m],\ \beta_{i}=\beta,\ i\in [m+1,a],\ \beta_{i}=0,\ i>a$.   Then using~(\ref{re1}), we have system of equations
  \begin{eqnarray}
  \label{a1}&&
  \alpha m+\beta (a-m)=1,\\
  && \alpha +\beta = 2\delta -2\frac{(1-2\delta)(k-1)}{n-2k} .
 \label{aa2} 
 \end{eqnarray}
  Assume that $\beta_{i}=\alpha, i\in [m],\ \beta_{i}=\beta,\ i\in [m+1,a],\ \beta_{i}=0,\ i>a$.  assume also that from at least one of equalities~(\ref{kk1}) follows equality
  \begin{equation}\label{ff1}
   \langle \beta ,x\rangle +\langle \beta ,\bar{x}\rangle=2\delta .
   \end{equation}
Hence there exist vectors $x ,\bar{ x}  \in {[n]\choose k}$, s.t. 
 \begin{eqnarray*}&& x_{j}=1,\ x_{a}=0,\ \bar{x}_j=0,\ \bar{x}_a =1,\\&&
 |\{x_{1},\ldots x_{m}\}| |= f_1,\   |\{x_{m+1},\ldots x_{a}\}| =f_2 ,\\&& |\{ \bar{x}_{1},\ldots ,\bar{x}_{m}\}| = \bar{f}_1,  |\{ \bar{x}_{m+1},\ldots ,\bar{x}_{a}\}| = \bar{f}_2 .
 \end{eqnarray*}
   We can rewrite equality~(\ref{ff1}) as follows 
   $$
   \alpha (f_1 +\bar{f}_1)+\beta (f_2+\bar{f}_2)=2\delta.
   $$
   Denote
   $$
   u=f_1+\bar{f}_1,\ \bar{u}=f_2+\bar{f}_2.
   $$
   Using equality~(\ref{aa2}) we have
   $$
   \alpha\left(2\omega -\frac{n-2k}{n-2}\right)+\beta \left(2\bar{u}-\frac{n-2k}{n-2}\right) =\frac{2(k-1)}{n-2}.
   $$
   From inequalities 
   $$
   2u-\frac{n-2k}{n-2} >u,\ 2\bar{u}-\frac{n-2k}{n-2} >\bar{u}.   $$
   we have
   $\delta <\frac{k-1}{n-2}$. From other side from~(\ref{aa2}) follows inequality
   $$
   \delta>\frac{k-1}{n-2}.
   $$
   This contradiction allows us to skip possibility of the case~(\ref{ff1}). 
   
    Hence we can assume that from~(\ref{kk1}) follows equalities
   $$
   \langle \beta_i ,x(r)\rangle = \langle \beta_i ,\bar{x}(r)\rangle .
   $$
   Summing these inequalities over $r$, we obtain relation
   $$\sum_{r=1}^{n-2\choose k-1}\langle \beta_i ,x(r)\rangle = \sum_{r=1}^{n-2\choose k-1}\langle \beta_i ,\bar{x}(r)\rangle 
   $$
   from which easily follows equality
   $$\alpha =\beta .$$
   Hence it is sufficient to consider the case $\alpha=\beta =\frac{1}{a}.$ 
    
We can rewrite inequalities in the definition of ${\cal A}$  in~(\ref{dd1}) as follows.

Let's 
$$
\beta_{i}\geq 0,\ \beta_i =0,\ i >k\ell-1;\ \sum_{i=1}^n \beta_i =1,\ \beta_{i}\geq \beta_{i+1} .
$$
Then the set
\begin{equation}
\label{e09}\left\{ x\in{[n]\choose k}:\ (\beta ,x) >\delta \right\}
\end{equation}
is compressed. 
It is easy to see, that arbitrary compressed set ${\cal A}$ is the intersection of finite number $N$ of compressed sets.

First we formulate the optimization problem:
\bigskip

{\bf Optimization Problem 2.}
\bigskip

Maximize over the choice of one-step functions $\beta $ and $\{\delta \}$ the function
$$
|{\cal A}|
$$
under the conditions that some  $s$  not $t$-intersection vectors ${x_t}$ does not belong to ${\cal A}$. The choice of $\{ x_t \}$ we will make belong.  

Solution and trust of this problem is literally the same as previous optimization problem, we only introduce the set $\{ x_t \}$.

Let  $a$ be the number of positive (equal) $\beta_{j}=1/a$. As in the previous section we leave only one restriction from the optimization problem 2. Restriction which is generated by the following $s$ elements:
\begin{eqnarray*}
x_1 &=& (\bar{I}_t,\hbox{shift}_0 (b),\ldots, \hbox{shift}_0 (b),z),\\
 x(i)&=& (\bar{I}_{t-1},0,\hbox{shift}_{i-1} (b),\ldots ,\hbox{shift}_{i-1} (b),z(i)),\\
\bar{I}_k &=&(1,\ldots ,1)\in \{0,1\}^k,\ z=(\bar{I}_{p},0,\ldots ,0)\in \{0,1\}^{n-t-fs};\\ && z(i)=(\hbox{shift}_{i-1} (\bar{I}_{p+1},0),\ldots ,0)\in\{0,1\}^{n-t-fs},\ i\in [2,s].
      \end{eqnarray*}
   Here $\hbox{shift}_i (w)$ is left cycle $i$-shifting of vector $w,\ b=(1,\ldots ,1,0)\in\{0,1\}^{s}$. We have $n=t+sf+g,$ where $g<s$ and $k=t+(s-1)f+p$. From other side from inequality
   $sk\leq (s-1)n+t$ it follows $p\leq g\frac{s-1}{s}$.
   
   We consider the restriction
   
 $$
 (\beta,x)>\delta
 $$
 which forbids this set of $n$-tuples.  It is easy to see, that  it is impossible that
  $a\in [t-1]$, because in this case should be $\delta \geq1$ and ${\cal A}=\emptyset$. If $a=t$, then we should choose $\delta =(t-1)/t .$
  This choice allows  $n$-tuple $x_1$ as a member of ${\cal A}$ and forbids other $x_i$. If
  $a=t+ms+r,\ r\in [s],\ m\leq f$, then we can choose some  $\delta \in ((a-s-1)/a, (a-s)/a) .$
  This choice of $\delta$  allows $x_1$ and forbids at least one  $x_i$. At last if $ (s-1)$ does not divide $(k-t)$, then for $a  \in [t+fs +r,n],$  the choice $\delta \in ((a -f-1)/a, (a-f)/a) $  forbids at least one $n$-tuple $x_i$. If $g>p+1$ then we can meet the situation when  $a\in [ t+fs+p+2,n]$, then choice  $\delta >k/a$ deliver set ${\cal A}=\emptyset $ because  allowing all $n$-tuples $x$ such that $|x\bigcap [a ]|=k$ does not warranty $s$-wise $t$-intersection, hence it is sufficient to consider the case $a \leq t+fs+p+1$. 
  
  Collecting these possibilities for the choice of of pairs $( a ,\delta )$ together we  see  that $(p<s)$
  \begin{equation}
  \label{ey7}
  N(s,n,k,t)\leq\max_{a=t+sp+r}\sum_{i\geq t+(s-1)p+r}  {a\choose i}{n-a\choose k-i} .
  \end{equation}
  It is sufficient to make the optimization in~(\ref{ey7}) only over $a$ such that
  $s|(a-t)$.
  Indeed it easily follows from the inequality
  $$
  \sum_{j\geq i-1}{a-1\choose j}{n-a+1\choose k-j}\geq\sum_{j\geq i}{a\choose j}{n-a\choose k-j}$$
  which can be proved by using the identity
  $$
  {k\choose m}={k-1\choose m}+{k-1\choose m-1}.
  $$
     
\end{document}